\documentclass[11pt,reqno]{amsart}
\usepackage{amssymb,latexsym}
\usepackage{graphicx}
\usepackage{url}		

\usepackage{color, url}
\usepackage[normalem]{ulem}
\usepackage{bm}

\calclayout

\makeatletter
\newcommand{\monthyear}[1]{%
  \def\@monthyear{\uppercase{#1}}}
\newcommand{\volnumber}[1]{%
  \def\@volnumber{\uppercase{#1}}}
\AtBeginDocument{%
\def\ps@plain{\ps@empty
  \def\@oddfoot{\@monthyear \hfil \thepage}%
  \def\@evenfoot{\thepage \hfil \@volnumber}}
\def\ps@firstpage{\ps@plain}
\def\ps@headings{\ps@empty
  \def\@evenhead{%
    \setTrue{runhead}%
    \def\thanks{\protect\thanks@warning}%
    \uppercase{The Fibonacci Quarterly}\hfil}%
  \def\@oddhead{%
    \setTrue{runhead}%
    \def\thanks{\protect\thanks@warning}%
    \hfill\uppercase{Hypergeometric Template}}%
  \let\@mkboth\markboth
  \def\@evenfoot{%
    \thepage \hfil \@volnumber}%
  \def\@oddfoot{%
    \@monthyear \hfil \thepage}%
  }%
\footskip=25pt
\pagestyle{headings}%
}
\makeatother

\newcommand{\G}{{\mathcal G}}
\newcommand{\Po}{{\mathcal P}}

\newcommand{\F}{{\mathcal{F}}}

\newcommand{\area}{{\texttt{area}}}
\newcommand{\sper}{{\texttt{sper}}}
\def\ver{{\textsf{ver}}}
\def\edg{{\textsf{edg}}}
\def\deg{{\textsf{deg}}}
\newcommand{\degv}{{\texttt{deg}}}
\newcommand{\Ham}{{\texttt{Ham}}}

\theoremstyle{plain}
\numberwithin{equation}{section}
\newtheorem{thm}{Theorem}[section]
\newtheorem{theorem}[thm]{Theorem}
\newtheorem{corollary}[thm]{Corollary}

\begin{document}
\monthyear{Month Year}
\volnumber{Volume, Number}
\setcounter{page}{1}

\title{Polyominoes and graphs
built from Fibonacci words}

\author{Sergey Kirgizov}
\address{LIB, Université de Bourgogne Franche-Comté \\
B.P. 47 870, 21078 Dijon Cedex, France}
\email{sergey.kirgizov@u-bourgogne.fr}

\author{Jos\'e L. Ram\'{\i}rez}
\address{Departamento de Matem\'aticas\\
 Universidad Nacional de Colombia,\\
   Bogot\'a, Colombia}
\email{jlramirezr@unal.edu.co}

\begin{abstract}

We introduce the $k$-bonacci polyominoes, a new family of polyominoes associated with the binary words avoiding  $k$ consecutive $1$'s, also called generalized $k$-bonacci words. 
The polyominoes are very entrancing objects, considered
in combinatorics and computer science. The study of polyominoes generates a rich source of combinatorial ideas.
In this paper we study some properties of $k$-bonacci polyominoes.
Specifically, we  determine their recursive structure 
and, using this structure, we enumerate them according
to their area, semiperimeter, and length
of the corresponding words.
We also introduce the $k$-bonacci graphs, then we obtain the generating functions for the total  number of  vertices and edges, the distribution of the degrees, and the total number of $k$-bonacci graphs that have a Hamiltonian cycle.  
\end{abstract}

\maketitle

The family of binary words avoiding a consecutive pattern is well-known in combinatorics and  computer science.  Many combinatorial statistics or parameters over these words can be studied with automata and grammars by means of the Chomsky-Sch\H{u}tzenberger methodology   \cite{DEL, flajolet, RAM}.
For instance, the language
of binary words avoiding $k$ consecutive $1$'s is a well-known example from one of Knuth's books~\cite[p. 286]{knuth3}. In this paper the 
words from this language are called \emph{generalized $k$-bonacci words}.    Vajnovszki~\cite{Vaj} studied these words in the context of exhaustive generation of Gray codes.  Recently, Bernini~\cite{Bernini}
considered some combinatorial properties  of these languages.  Baril et al.~\cite{Baril} give a bijection between $k$-bonacci words and the $q$-decreasing words
(for $q = k-1$), they also provide an efficient exhaustive generating algorithm for $q$-decreasing words in lexicographic order.  

Another direction in the study of the Fibonacci words is in Graph Theory. The $n$-length binary words that avoid two consecutive ones are the vertices of the Fibonacci cube~\cite{Hsu}. Two Fibonacci words  of the same length  are adjacent in the graph if its Hamming distance is equal to one, that is, they differ in exactly one symbol. The Fibonacci cube is a subgraph of the $n$-dimensional hypercube. The Fibonacci cube has been extensively studied 
in recent years. See~\cite{Klav2} for a survey.

In this paper, we are interested in the study of a new family of polyominoes and graphs associated with the generalized $k$-bonacci words. Let $\F_{n,k}$ denote the set of $n$-length binary words avoiding  $k$ consecutive $1$'s, and $\F_{k}=\cup_{n\geq 1}\F_{n,k}$. The set $\F_{k}$ corresponds to the set of generalized $k$-bonacci words.  The elements of $\F_{n,2}$ are called \emph{Fibonacci words}.  
  For example,
$$\F_{3,2}=\{000, 001, 010, 100, 101\} \quad \text{and} \quad  \F_{3,3}=\{000, 001, 010, 011, 100, 101, 110\}.$$

\begin{figure}[ht]
\centering
\includegraphics[scale=1]{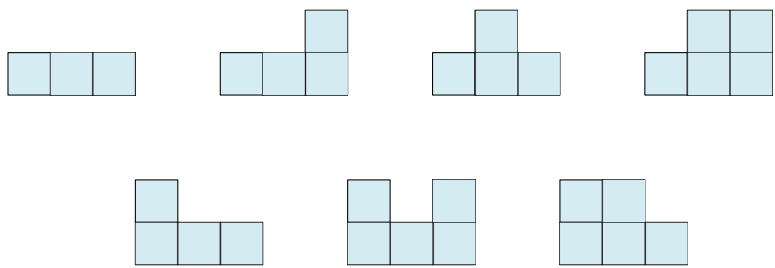}
\caption{The Fibonacci polyominoes associated with the words in $\F_{3,3}$.}
\label{fig1}
\end{figure}

The set $\F_{n,k}$ is enumerated   by the  generalized Fibonacci numbers $F_{n+2,k}$. This sequence is defined by $F_{n,k}=\sum_{i=1}^kF_{n-i,k}$ for $n\geq 2$, with $F_{1,k}=1$ and $F_{n,k}=0$ for all $n\leq 0$. Given a word $w=w_1\cdots w_n\in \F_{n,k}$, its associated polyomino, called \emph{$k$-bonacci polyomino}, is a bargraph of  $n$ columns, such that the $i$-th column has $w_i+1$ unit cells for $1\leq i \leq n$. For example, Figure~\ref{fig1} shows the polyominoes associated with the  $3$-bonacci words of length $3$.  Let $\Po_{n,k}$ denote the set of $k$-bonacci polyominoes with $n$ columns. The elements of $\Po_{n,2}$ are called \emph{Fibonacci polyominoes}. 

 The polyominoes provide a rich source of combinatorial ideas. For example, polyominoes play an important role in the combinatorics on words because they can be encoded by words, and then the problem of deciding if a given polyomino tiles the plane  by a translation  reduces to finding a special factorization of the word~\cite{NIV}.  In particular, polyominoes associated with some special prefixes of the infinite  Fibonacci word  tile the plane by translations~\cite{BLO2, RAMRUB}.  The polyominoes  have also been studied in connection with other  discrete structures such as permutations, set partitions,  compositions,  among others (see for example~\cite{BLE3, CallManRam, Book1, ManSha} and references contained therein).

\begin{figure}[ht]
\centering
\includegraphics[scale=1]{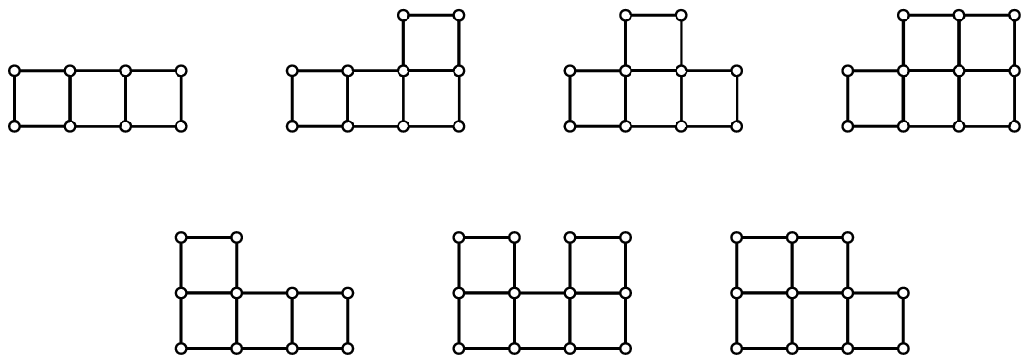}
\caption{The $3$-bonacci graphs  in $\G_{3,3}$.}
\label{fig2}
\end{figure}

 Any $k$-bonacci polyomino can be  regarded as a 
graph, called  \emph{$k$-bonacci graph},  considering  the cell sides as edges and cell corners
as vertices. For 
example, Figure~\ref{fig2} shows the graphs  
associated with the  $3$-bonacci polyominoes 
with $3$ columns.  Let $\G_{n,k}$ denote the set
of $k$-bonacci graphs associated with the 
polyominoes of $\Po_{n,k}$.  Note that in 
$\G_{n,k}$ there are isomorphic graphs, in 
particular, graphs corresponding to words 00
and 1 are isomorphic. And if $w_1$ and $w_2$ are two different $k$-bonacci words of length $n\geq 2$,  such that
$w_1^R=w_2$ 
(where $w^R$ denotes the 
reverse of the word $w$), then the graphs induced by these 
words are isomorphic.
The $k$-bonacci 
graphs are particular examples of chemical 
graphs, that is, graphs  with all vertices of 
degree at most four.  

In this paper, we obtain several enumerative results  of the $k$-bonacci polyominoes    including the area and  semi-perimeter.  For the $k$-bonacci graphs we study the number of vertices and edges, the degree sequence polynomial,  the average degree of a vertex, and  the number of $k$-bonacci graphs that have at least one Hamiltonian cycle.

We obtain these results by using the recursive construction
of the $k$-bonacci words  and  generating functions in several variables. The generating functions  
have been successfully used to study several 
statistics over Fibonacci-runs graphs and the 
restricted Fibonacci words.
See for example~\cite{OME1, OME2, OME3}. 
 
\section{Decomposition of the $k$-Bonacci Polyominoes}

\begin{figure}[ht]
\centering
\includegraphics[scale=1]{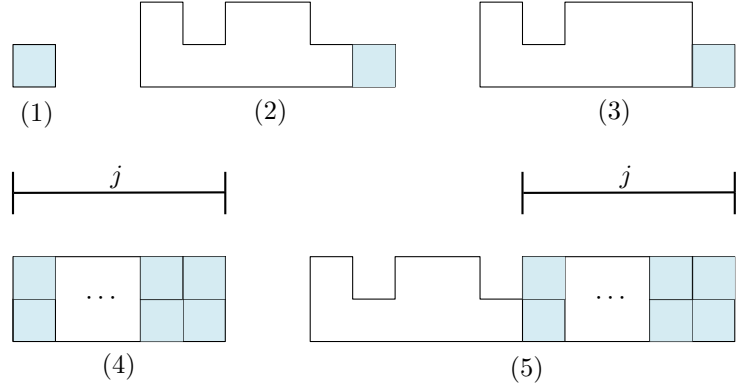}
\caption{Decomposition of a $k$-bonacci polyomino.}
\label{deco}
\end{figure}

Let $\Po_{n,k}^{(1)}$ and  $\Po_{n,k}^{(2)}$ denote the sets of $k$-bonacci polyominoes with $n$ columns, whose last column has height $1$ and 2, respectively. It is clear that $\Po_{n,k}=\Po_{n,k}^{(1)} \cup \Po_{n,k}^{(2)}$. If $P\in \Po_{n,k}^{(1)}$, then $P$ is a unit square or $P=P'\square$, where $P' \in \Po_{n-1,k}=\Po_{n-1,k}^{(1)} \cup \Po_{n-1,k}^{(2)}$ $(n>1)$.  If $P\in \Po_{n,k}^{(2)}$, then $P=P'C_j$, where $P' \in \Po_{n-j,k}^{(1)}$ (possibly empty) and  $C_j$ is a concatenation of $j$ columns of height 2, for $1\leq j < k$.  Figure~\ref{deco} illustrates the above decomposition.

\subsection{Area and Semiperimeter}

In this section we study the distribution of the area and semiperimeter    in $\Po_{n,k}$. 
Let $P$ be a $k$-bonacci polyomino.  We denote  by $\area(P)$  the number of cells of $P$ and by
$\sper(P)$ the semiperimeter, that is the half of  the  perimeter of $P$ (since  the perimeter of a $k$-bonacci polyomino is always an even number).  

\begin{figure}[ht]
\centering
\includegraphics[scale=1]{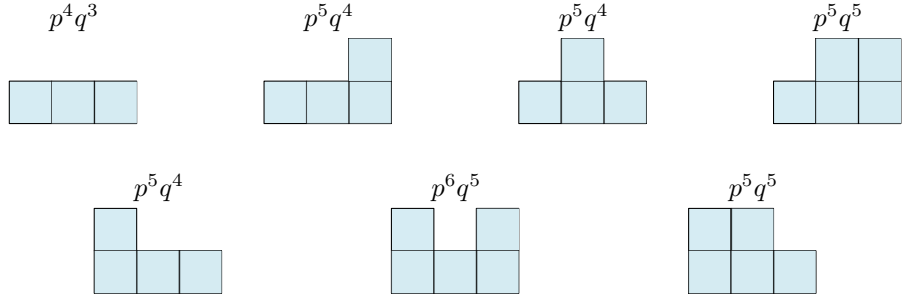}
\caption{Weights for polyominoes in $\Po_{3,3}$.}
\label{grafo2}
\end{figure}

 Define the generating function
$$F_k(x,p,q):=\sum_{n\geq 1} x^n \sum_{P\in \Po_{n,k}}p^{\sper(P)}q^{\area(P)}, \quad k\geq 2,$$
where $x$ marks the length of the corresponding $k$-bonacci word, i.e., the number of cells in the bottom row of a polyomino.
Analogously, we introduce the generating functions
$$F_{k,j}(x,p,q):=\sum_{n\geq 1} x^n \sum_{P\in \Po_{n,k}^{(j)}}p^{\sper(w)}q^{\area(P)},  \quad \text{for} \quad j=1, \,  2.$$
It is clear that 
\begin{align}\label{GF}
F_k(x,p,q)=F_{k,1}(x,p,q)+F_{k,2}(x,p,q).
\end{align}
In Theorem~\ref{teoareaper1} we give a rational expression for the generating function $F_k(x,p,q)$. 
\begin{theorem}\label{teoareaper1}
The generating function for $k$-bonacci polyominoes with respect to
the length of the corresponding word,
semiperimeter and area (marked respectively by $x$, $p$ and $q$)
is
$$F_k(x,p,q)=\frac{p^2 ((q + p q^2 ) x - (p q^3 - p^2 q^3) x^2 - p^k q^{2 k} x^k - 
   p^{k+1} q^{2 k+1} x^{k+1})}{1 - (p q + p q^2) x + (p^2 q^3 - p^3 q^3) x^2 + 
 p^{k+2} q^{2 k+1} x^{k+1}}.$$
\end{theorem}
\begin{proof}
From the decomposition given in Figure~\ref{deco} we have the functional equations
\begin{align*}
F_{k,1}(x,p,q)&=\underbrace{p^2qx}_{(1)} + \underbrace{pqx(F_{k,1}(x,p,q)+F_{k,2}(x,p,q))}_{(2)+(3)} \\
F_{k,2}(x,p,q)&=\underbrace{\sum_{j=1}^{k-1}p^{i+2}q^{2i}x^j}_{(4)} +\underbrace{\left(\sum_{j=1}^{k-1}p^{i+1}q^{2i}x^j\right) F_{k,1}(x,p,q)}_{(5)}\\
&=\frac{p^2(pq^2x-p^kq^{2k}x^k)}{1-pq^2x} + \frac{p(pq^2x-p^kq^{2k}x^k)}{1-pq^2x}F_{k,1}(x,p,q).
\end{align*}
To make understanding easier, cases (1) to (5) in Figure~\ref{deco} are indicated below their corresponding terms.
Solving the above system of equation and from \eqref{GF} we obtain the desired result.
\end{proof}
For example, the series expansion of the generating function $F_{3}(x,p,q)$ begins with
\begin{multline*}
(p^2 q + p^3 q^2) x + (p^3 q^2 + 2 p^4 q^3 + 
    p^4 q^4) x^2 + \bm{(p^4 q^3 + 3 p^5 q^4 + 2 p^5 q^5 + 
    p^6 q^5) x^3} \\
    + (p^5 q^4 + 4 p^6 q^5 + 3 p^6 q^6 + 3 p^7 q^6 +  2 p^7 q^7) x^4 + \cdots
    \end{multline*}
Figure~\ref{grafo2} shows the weights of the  $3$-bonacci polyominoes   corresponding to the bold coefficient in the above series.

\begin{corollary}\label{coroarea}
The generating functions for the total area and semi-perimeter within all of the members in $\Po_{n, k}$ are 
$$A_k(x)=\frac{3 x - 2 k x^k - 2 (2 - k) x^{k+1} + x^{2 k+1}}{(1 - 2 x + x^{k+1})^2}$$
and 
$$P_k(x)=\frac{5 x - 5 x^2 - (2 + k) x^k - (1 - k) x^{k+1} + 4 x^{k+2} - x^{2k+2}}{(1 - 2 x + x^{k+1})^2},$$
respectively. 
\end{corollary}
\begin{proof}
From the definition of  $F_k(x,p,q)$ we have the equalities
$$A_k(x)=\left.\frac{\partial F_k(x,1,q)}{\partial q}\right|_{q=1} \quad\text{and}\quad P_k(x)=\left.\frac{\partial F_k(x,p,1)}{\partial p}\right|_{p=1}.$$
Therefore, from Theorem~\ref{teoareaper1} we obtain the desired results.
\end{proof}
In particular, for $k = 2, 3,$ and $4$ we obtain the following generating functions:
\begin{align*}
A_2(x)&=\frac{x (3 + 2 x + x^2)}{(1 -  x -x^2)^2},\\
 A_3(x)&=\frac{x (3 + 6 x + 3 x^2 + 2 x^3 + x^4)}{(1 -  x -x^2-x^3)^2},\\
 A_4(x)&=\frac{x (3 + 6 x + 9 x^2 + 4 x^3 + 3 x^4 + 2 x^5 + x^6)}{(1 -  x -x^2-x^3-x^4)^2},\\
P_2(x)&=\frac{x (5 + x - 2 x^2 - x^3)}{(1 -  x -x^2)^2},\\
 P_3(x)&=\frac{x (5 + 5 x - 3 x^3 - 2 x^4 - x^5)}{(1 -  x -x^2-x^3)^2},\\
 P_4(x)&=\frac{x (5 + 5 x + 5 x^2 - x^3 - 4 x^4 - 3 x^5 - 2 x^6 - x^7)}{(1 -  x -x^2-x^3-x^4)^2}.
   \end{align*}

We will apply the same technique  to obtain other corollaries  throughout the article.

\subsection{Area and Semiperimeter of the Fibonacci polyominoes}

For the Fibonacci polyominoes ($k=2$) we can give some additional results.  Let $t_n(p,q)$ denote the $n$-th coefficient of  $F_2(x,p,q)$, that is, 
$$t_n(p,q):= \sum_{P\in \Po_{n,2}}p^{\sper(P)}q^{\area(P)}.$$

\begin{theorem}\label{fibopoly}
For all $n\geq 3$ we have
$$t_n(p,q)=pqt_{n-1}(p,q) + p^3q^3t_{n-2}(p,q),$$
with the initial values $t_1(p,q)=p^2 q + p^3 q^2$ and $t_2(p,q)=p^3 q^2 + 2 p^4 q^3$. Moreover,  for all $n\geq 1$ we have the combinatorial formula 
\begin{align}\label{closedf2b}
t_n(p,q)=\sum_{i=0}^{\lfloor \frac{n+1}{2} \rfloor} \binom{n+1-i}{i}p^{n + i + 1} q^{n + i}.
\end{align}
\end{theorem}
\begin{proof}
Let $P$ be a Fibonacci polyomino with $n$ columns ($n\geq 3$).  If the last column has height 1, then the number of this kind of polyominoes is given by  $qpt_{n-1}(p,q)$. On the other hand, if the last column has height 2, then the previous column has height 1, and these polyominoes are counted by the polynomial  $p^3q^3t_{n-2}(p,q)$.  Hence the total number of Fibonacci polyominoes is given by $pqt_{n-1}(p,q) + p^3q^3t_{n-2}(p,q)$. Finally, we can use mathematical induction and the recurrence relation to prove the combinatorial identity.  This formula can be also proved by means of Zeilberger's creative telescoping method~\cite{ABZ}. Denote the summand on the right side   
of the equality in \eqref{closedf2b} by $F(n,i)$, that is
$$F(n,i):= \binom{n+1-i}{i}p^{n + i + 1} q^{n + i}.$$
By the Zeilberger algorithm, $F(n,i)$ satisfies the relation
\begin{align}\label{rel1}
F(n+2,i)-pqF(n+1,i)-p^3q^3F(n,i)=G(n,i+1)-G(n,i),
\end{align}
with the certificate 
$$R(n,i)=-\frac{i (2 - i + n) p^2 q^2}{(2 - 2 i + n) (3 - 2 i + n)}.$$
That is, $R(n,i)=G(n,i)/F(n,i)$ is a rational function in both variables.
If $f(n)$ denotes the right sum in the equality \eqref{closedf2b}, then summing both sides of \eqref{rel1} with respect to $i$ yields $f(n+2)-pqf(n+1)-p^3q^3f(n)=0$.  The sequences $f(n)$ and $t_n(p,q)$ satisfy the same recurrence relation and  have the same initial values, therefore these sequences coincides for all positive integers $n$.

\end{proof}

\begin{corollary}
The  total area for all Fibonacci polyominoes in $\Po_{n,2}$ is 
$$\sum_{i=0}^{\lfloor \frac{n+1}{2} \rfloor} \binom{n+1-i}{i}(n + i)=\frac{1}{5}(6nF_{n+2}+(n+2)F_n),$$
where $F_n$ is the $n$-th Fibonacci number.
\end{corollary}
\begin{proof}
From Theorem~\ref{fibopoly} the total area of the Fibonacci polyominoes is given by
$$\left.\frac{\partial t_n(1,q)}{\partial q}\right|_{q=1} =\sum_{i=0}^{\lfloor \frac{n+1}{2} \rfloor} \binom{n+1-i}{i}(n + i).$$
For the expression in terms of Fibonacci numbers we use the generating function of the total area
$$A_2(x)=\frac{x (3 + 2 x + x^2)}{(1 -  x -x^2)^2}=(3x+2x^2+x^3)\sum_{n\geq 0}c(n+2)x^n,$$
where $c(n)$ is the convolution of the Fibonacci numbers, that is 
$$c(n)=\sum_{i=0}^nF_iF_{n-i}=\frac{1}{5}\left((n-1)F_n + 2nF_{n-1}\right).$$
In the last equality we use the identity (32.13) given in~\cite{Koshy}. Therefore,  
$$[x^n]A_2(x)=3c(n+1)+2c(n)+c(n-1)=\frac{1}{5}(6nF_{n+2}+(n+2)F_n), \quad  n \geq 1.$$
\end{proof}

The number of Fibonacci polyominoes of area $n$ is related to the 
Narayana's cows sequence ~\cite{Allou}. The 
Narayana's cows sequence $b_n$ is defined by the recurrence relation  $b_{n}=b_{n-1}+b_{n-3}$ for $n\geq3$, with the initial values  $b_0=0, b_1=1$, and  $b_2=1$. The Narayana's cows sequence can be calculated with the formula 
$$b_n=\sum_{i=0}^{\lfloor (n-1)/3 \rfloor} \binom{n-2i-1}{i}.$$

\begin{theorem}
The number of Fibonacci polyominoes of area $n$ is equal to 
 the number of Narayana's cows $b_{n+1}$. 
\end{theorem}
\begin{proof}
From Theorem~\ref{teoareaper1} the generating function of the number of Fibonacci polyominoes of a fixed area is 
$$F_2(1,1,q)=\frac{q (1 + q + q^2)}{1-q-q^3}.$$
On the other hand, the generating function of the Narayana sequence is  $N(x):=\sum_{i\geq 0}b_ix^i=1/(1-x-x^3)$. Then 
$N(q)=1+q+qF_2(1,1,q)$. By comparing
the $n$-th coefficient of the generating functions $N(q)$ and $F_2(1,1,q)$ we obtain  the desired result. 
\end{proof}

\begin{figure}[ht]
\centering
\includegraphics[scale=1]{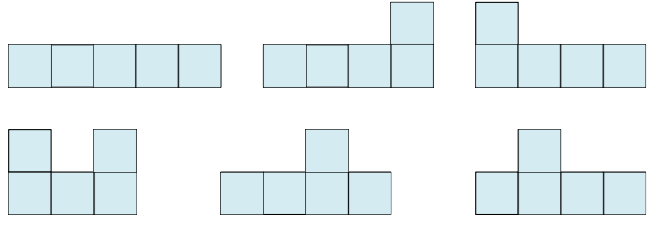}
\caption{Fibonacci polyominoes of area 5.}
\label{grafo5}
\end{figure}

For example,  Figure~\ref{grafo5} shows the Fibonacci polyominoes   of area 5, that is, $b_{6}=6$ polyominoes. 

\section{Number of Vertices and Edges}

The goal of this section is to enumerate the number of vertices and edges of the $k$-bonacci graphs. 
Let $G$ be a $k$-bonacci graph.  We denote  by $\ver(G)$ and $\edg(G)$ the number of vertices and edges  of the graph $G$.  Let $\G_{n,k}^{(1)}$ and  $\G_{n,k}^{(2)}$ denote the set of $k$-bonacci graphs associated with the polyominoes in  $\Po_{n,k}^{(1)}$ and  $\Po_{n,k}^{(2)}$, respectively.
   Define the generating function
$$G_k(x,p,q):=\sum_{n\geq 1} x^n \sum_{G\in \G_{n,k}}p^{\edg(G)}q^{\ver(G)}, \quad k\geq 2.$$
Similarly, we have the generating functions
$$G_{k,j}(x,p,q):=\sum_{n\geq 1} x^n \sum_{G\in \G_{n,k}^{(j)}}p^{\edg(G)}q^{\ver(G)},  \quad \text{for} \quad j=1, 2.$$
In these trivariate generating functions
the variable $x$ marks the length of the corresponding $k$-bonacci word, i.e., the number of vertices in the bottom row of a graph minus one.

\begin{figure}[ht]
\centering
\includegraphics[scale=1]{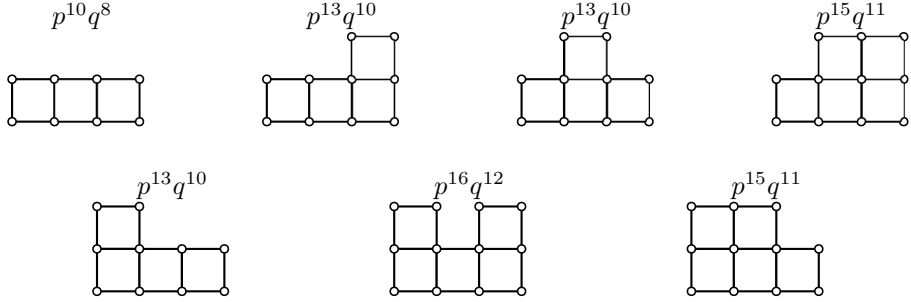}
\caption{Weights for the graphs in $\G_{3,3}$.}
\label{grafo3}
\end{figure}

\begin{theorem}
The generating function for $k$-bonacci graphs with respect to 
the length of the corresponding word, 
the number of edges and vertices (marked respectively by $x$, $p$ and $q$)
is
$$G_k(x,p,q)=\frac{p^2 q^3((p^2 q + p^5 q^3) x - (p^7 q^4 - p^8 q^5) x^2 - p^{5k} q^{3 k} x^k -
   p^{5 k+3} q^{3 k+2} x^{k+1})}{1 - (p^3 q^2 + p^5 q^3) x + (p^8 q^5 - p^9 q^6) x^2 + 
 p^{5 k+4} q^{3 k+3} x^{k+1}}.$$

\end{theorem}
\begin{proof}
From the decomposition given in Figure~\ref{deco} we have the functional equations
\begin{align*}
G_{k,1}(x,p,q)&=\underbrace{p^4q^4x}_{(1)} + \underbrace{p^3q^2x(G_{k,1}(x,p,q)+G_{k,2}(x,p,q))}_{(2)+(3)} \\
G_{k,2}(x,p,q)&=\underbrace{\sum_{j=1}^{k-1}p^{3 j + 2 j + 2}q^{3 (j + 1)}x^j}_{(4)} +\underbrace{\left(\sum_{j=1}^{k-1}p^{3 j + 2 j + 1}q^{3 j + 1}x^j\right) G_{k,1}(x,p,q)}_{(5)}\\
&=\frac{p^2 q^3 (p^5 q^3 x - p^{5 k}q^{3 k} x^k)}{1 - p^5 q^3 x} + \frac{p q (p^5 q^3 x - p^{5 k}q^{3 k} x^k)}{1 - p^5 q^3 x}G_{k,1}(x,p,q).
\end{align*}
Solving the above system of equation and from the equality $G_k(x,p,q)=G_{k,1}(x,p,q)+G_{k,2}(x,p,q)$ we obtain the desired result.
\end{proof}
For example, the series expansion of the generating function $G_{3}(x,p,q)$ begins with
\begin{multline*}
x \left(p^7 q^6+p^4 q^4\right)+x^2 \left(p^{12} q^9+2 p^{10}
   q^8+p^7 q^6\right)
   + \bm{x^3 \left(p^{16} q^{12}+2 p^{15} q^{11}+3
   p^{13} q^{10}+p^{10} q^8\right)}\\
   +x^4 \left(2 p^{21} q^{15}+3 p^{19} q^{14}+3 p^{18} q^{13}+4 p^{16} q^{12}+p^{13}
   q^{10}\right)+ \cdots
    \end{multline*}
Figure~\ref{grafo3} shows the weights of the  $3$-bonacci graphs  corresponding to the bold coefficient in the above series.

\begin{corollary}
The generating functions for the total number of vertices and edges  within all of the members in $\G_{n, k}$ are 
$$V_k(x)=\frac{10 x - 9 x^2 - 3 (1 + k) x^k - (4 - 3 k) x^{k+1} + 8 x^{k+2} -  2 x^{2k+2}}{(1 - 2 x + x^{k+1})^2}$$
and
$$E_k(x)=\frac{11 x - 5 x^2 - (2 + 5 k) x^k - (9 - 5 k) x^{k+1} + 4 x^{k+2} + 2 x^{2 k+1} - x^{2 k+2}}{(1 - 2 x + x^{k+1})^2},$$
respectively.
\end{corollary}

In particular, for $k = 2, 3$, and $4$ we obtain the following generating functions:
\begin{align*}
V_2(x)&=\frac{2 x (5 + x - 2 x^2 - x^3)}{(1 -  x -x^2)^2},\\
 V_3(x)&=\frac{x (10 + 11 x - 6 x^3 - 4 x^4 - 2 x^5)}{(1 -  x -x^2-x^3)^2},\\
 V_4(x)&=\frac{x (10 + 11 x + 12 x^2 - 2 x^3 - 8 x^4 - 6 x^5 - 4 x^6 - 2 x^7)}{(1 -  x -x^2-x^3-x^4)^2}, \\
E_2(x)&=\frac{x (11 + 5 x - x^3)}{(1 -  x -x^2)^2},\\
 E_3(x)&=\frac{x (11 + 17 x + 6 x^2 + x^3 - x^5)}{(1 -  x -x^2-x^3)^2},\\
 E_4(x)&=\frac{x (11 + 17 x + 23 x^2 + 7 x^3 + 2 x^4 + x^5 - x^7)}{(1 -  x -x^2-x^3-x^4)^2}.
   \end{align*}

Let $v_n(p,q)$ denote the $n$-th coefficient of  $G_2(x,p,q)$, that is, 
$$v_n(p,q):= \sum_{G\in \G_{n,2}}p^{\edg(G)}q^{\ver(G)}.$$

\begin{theorem}\label{fibopoly2}
For all $n\geq 3$ we have
$$v_n(p,q)=p^3q^2v_{n-1}(p,q) + p^9q^6v_{n-2}(p,q),$$
with the initial values $v_1(p,q)=p^4 q^4 + p^7 q^6$ and $v_2(p,q)=p^7 q^6 + 2 p^10 q^8$. Moreover,  for all $n\geq 1$ we have the combinatorial formula 
\begin{align*}
v_n(p,q)=\sum_{i=0}^{\lfloor \frac{n+1}{2} \rfloor} \binom{n+1-i}{i}p^{3 n + 1 + 3 i} q^{2 n + 2 + 2 i}.
\end{align*}
\end{theorem}

\section{Degree Sequences}

In this section we are interested in the degree sequence of a $k$-bonacci graph. Degree sequences have been well studied for Fibonacci cubes~\cite{Klav5} and Fibonacci-run graphs~\cite{OME3}. 
Recall that the degree of a vertex of a graph is the number of edges that are incident to the vertex.  A $k$-bonacci graph can have only vertices of degree two, three, or four.  Let $G$ be a $k$-bonacci graph. We denote  by $\degv_i(G)$ the number of vertices of degree $i$ in the graph  $G$. 

We are interested in the  generating function
$$D_k(x, q_2,q_3,q_4):=\sum_{n\geq1}x^n\sum_{G\in \G_{n,k}}q_2^{\degv_2(G)}q_3^{\degv_3(G)}q_4^{\degv_4(G)}, \quad k\geq 2,$$
 where $x$ marks the length of the corresponding $k$-bonacci word, i.e., the number of vertices in the bottom row of a graph minus one.
Analogously, we introduce  the generating functions
$$D_{k,j}(x, q_2,q_3,q_4):=\sum_{n\geq 1} x^n \sum_{G\in \G_{n,k}^{(j)}}q_2^{\degv_2(G)}q_3^{\degv_3(G)}q_4^{\degv_4(G)}, \quad j=1,  2.$$

\begin{theorem}\label{degreegfg}
For all $k\geq 2$, the generating function $D_k(x,q_2,q_3,q_4)$ is given by 
$$\frac{q_2^4 \left( (q_3^2 q_4 +q_4) x -(q_3^2
   q_4^2 - 2 q_2 q_3^2 q_4^2 +q_3^4 q_4) x^2 - q_3^{2 k} q_4^k x^k + (q_3^{2 k+2} q_4^k -2 q_2 q_3^{2 k}
   q_4^{k+1}) x^{k+1}\right)}{q_4(1 -(q_3^2 + q_4 q_3^2) x + (q_4 q_3^4 - q_2^2 q_4^2 q_3^2) x^2 + q_2^2 q_4^{k+1} q_3^{2 k}
   x^{k+1})}.$$
\end{theorem}
\begin{proof}
Let $G$ be a Fibonacci graph in $\G_{n,k}$.  If $n=1$, then $G$ contributes to the generating function the term $q_2^4x$. See Figure~\ref{deco2} case (1).  If $n>1$ and $G\in \G_{n,k}^{(1)}$, then this case contributes to the generating function the terms $$q_3^2xD_{k,1}(x,q_2,q_3,q_4) \quad \text{and} \quad q_2q_4xD_{k,2}(x,q_2,q_3,q_4),$$  as seen in Figure~\ref{deco2} cases (2) and (3). Notice that in the case (2) we have only two new vertices of degree 3.  In the case (3) we have a new vertex of degree 2 and another of degree 4. The red vertices denote the vertices of degree 2, the blue vertices denote the vertices of degree 3, and the green vertices denote the vertices of degree 4. For $n>1$ and $G\in \G_{n,k}^{(2)}$, this case contributes the terms   (see Figure~\ref{deco2} cases (4) and (5))
\begin{align*}
\underbrace{q_2^4\sum_{j=1}^{k-1}q_3^{2 (j - 1) + 2}q_4^{j - 1}x^j}_{(4)} +\underbrace{\left(q_2\sum_{j=1}^{k-1}q_3^{2 (j - 1) + 2}q_4^{j}x^j\right) D_{k,1}(x,q_2,q_3,q_4)}_{(5)}.
\end{align*}

\begin{figure}[ht]
\centering
\includegraphics[scale=1]{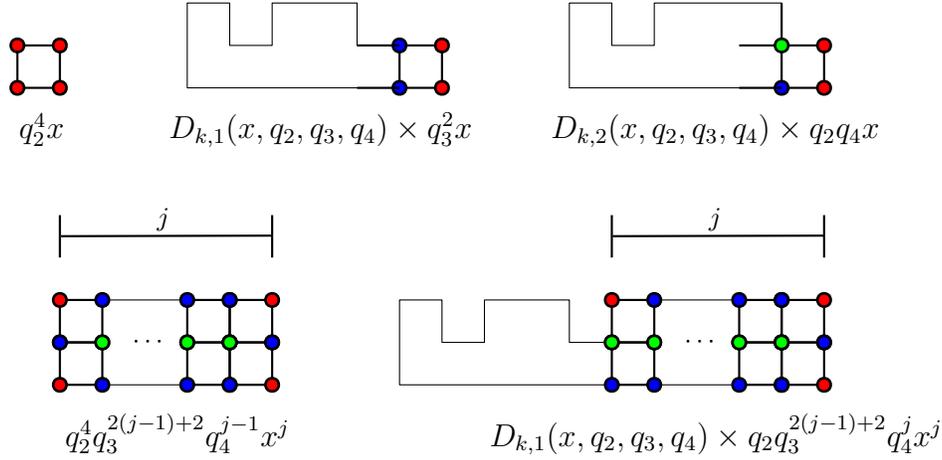}
\caption{Decomposition of a $k$-bonacci graph.}
\label{deco2}
\end{figure}
Therefore, we have the functional equations
\begin{align*}
D_{k,1}(x,q_2,q_3,q_4)&=q_2^4x+ q_3^2xD_{k,1}(x,q_2,q_3,q_4) +  q_2q_4xD_{k,2}(x,q_2,q_3,q_4)\\
D_{k,2}(x,q_2,q_3,q_4)&=\frac{q_2^4 (q_3^2q_4 x - q_3^{2k}q_4^k x^k)}{q_4(1-q_3^2q_4x)} +\frac{q_2 (q_3^2q_4 x - q_3^{2k}q_4^k x^k)}{1-q_3^2q_4x}  D_{k,1}(x,q_2,q_3,q_4).
\end{align*}

Solving the above system of equation we obtain the desired result.
\end{proof}
For example, the series expansion of the generating function $D_3(x,q_2,q_3,q_4)$ begins with
\begin{multline*}
\left(q_2^4q_3^2+q_2^4\right) x+\left(2 q_2^5 q_3^2 q_4 + q_2^4 q_3^2
   + q_2^4 q_3^4 q_4 \right) x^2+\\
   \bm{\left(q_2^6 q_3^4 q_4^2 + q_2^6 q_3^2
   q_4^2 +2 q_2^5 q_3^4 q_4^2 + 2 q_2^5 q_3^4 q_4  + q_3^4 q_2^4\right) x^3} + \cdots
    \end{multline*}
Figure~\ref{grafo4} shows the weights of the  $3$-bonacci graphs  corresponding to the bold coefficient in the above series. 
\begin{figure}[ht]
\centering
\includegraphics[scale=1]{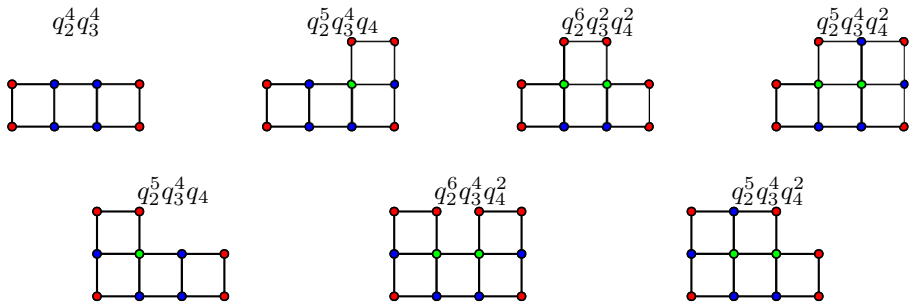}
\caption{Weights for the graphs in $\G_{3,3}$.}
\label{grafo4}
\end{figure}

\begin{corollary}
The generating function for the total number of vertices of degree 2 within all of the members in $\G_{n, k}$ is
$$D_{k}^{(2)}(x):=\frac{2 (4 x - 7 x^2 - 2 x^k + x^{k+1} + 7 x^{k+2} - x^{2 k+1} - 2 x^{2 k+2})}{(1 - 2 x + x^{k+1})^2}.$$
\end{corollary}

In particular, for $k = 2, 3$, and $4$ we obtain the following generating functions:
\begin{align*}
D_{2}^{(2)}(x)&=\frac{2 x (4 - x - 5 x^2 - 2 x^3)}{(1 -  x -x^2)^2},\\
D_{3}^{(2)}(x)&=\frac{2 x (4 + x - 4 x^2 - 8 x^3 - 5 x^4 - 2 x^5)}{(1 -  x -x^2-x^3)^2},\\
D_{4}^{(2)}(x)&=\frac{2 x (4 - 7 x - 2 x^3 + x^4 + 7 x^5 - x^8 - 2 x^9)}{(1 - 2 x + x^5)^2}.
   \end{align*}
\begin{corollary}
The generating function for the total number of vertices of degree 3 within all of the members in $\G_{n, k}$ is
$$D_{k}^{(3)}(x):=\frac{2 (x + x^2 - k x^k - (1 - k) x^{k+1} - 2 x^{k+2} + x^{2k+2})}{(1 - 2 x + x^{k+1})^2}.$$
\end{corollary}
In particular, for $k = 2, 3$ and $4$ we obtain the following generating functions:
\begin{align*}
D_{2}^{(3)}(x)&=\frac{2 x (1 + x + 2 x^2 + x^3)}{(1 -  x -x^2)^2},\\
D_{3}^{(3)}(x)&=\frac{2 (x + x^2 - 3 x^3 + 2 x^4 - 2 x^5 + x^8)}{(1 - 2 x + x^4 )^2},\\
D_{4}^{(3)}(x)&=\frac{2(x + x^2 - 4 x^4 + 3 x^5 - 2 x^6 + x^{10})}{(1 - 2 x + x^5)^2}.
   \end{align*}
   
  \begin{corollary}
The generating function for the total number of vertices of degree 4 within all of the members in $\G_{n, k}$ is
$$D_{k}^{(4)}(x):=\frac{3 x^2 + (1-k) x^k - (4 - k) x^{k+1} - 2 x^{k+2} + 2x^{2k+1}}{(1 - 2 x + 2 x^{k+1})^2}.$$
\end{corollary}

In particular, for $k = 2, 3$ and $4$ we obtain the following generating functions:
\begin{align*}
D_{2}^{(4)}(x)&=\frac{2 x^2 (1 + x)}{(1 -  x -x^2)^2},\\
D_{3}^{(4)}(x)&=\frac{x^2 (3 + 4 x + 4 x^2 + 2 x^3)}{(1 - x - x^2 - x^3 )^2},\\
D_{4}^{(4)}(x)&=\frac{3 x^2 - 3 x^4 - 2 x^6 + 2 x^9}{(1 - 2 x + x^5)^2}.
\end{align*}
   
\subsection{The degree sequence of the Fibonacci graph}
Finally, we consider the particular case of Fibonacci graphs ($k=2$).  Let $d_{n,2}(q)$ denote the $n$-th coefficient of  $D_2(x, q,1,1)$, that is, 
$$d_{n,2}(q):= \sum_{G\in \G_{n,2}}q^{\deg_2(G)}.$$
Similarly, we define the sequences
$$d_{n,3}(q):= \sum_{G\in \G_{n,2}}q^{\deg_3(G)} \quad \text{and} \quad d_{n,4}(q):= \sum_{G\in \G_{n,2}}q^{\deg_4(G)}.$$
\begin{theorem}\label{closeformula}
For all $n\geq 3$ we have
$$d_{n,2}(q)=d_{n-1,2}(q) + q^2d_{n-2,2}(q),$$
with the initial values $d_{1,2}(q)=2q^4$ and $d_{2,2}(q)=q^4 + 2 q^4$. Moreover,  for all $n\geq 1$ we have the combinatorial formula 
\begin{align}\label{vertdeg2}
d_{n,2}(q)=\sum_{i=1}^{n} \left(\binom{n-1-\lfloor \frac{i}{2}\rfloor}{\lfloor\frac{i-1}{2}\rfloor}+  \binom{n-2-\lfloor \frac{i-1}{2}\rfloor}{\lfloor\frac{i-2}{2}\rfloor}\right)q^{i+3}.
\end{align}
\end{theorem}
\begin{proof}
Define the generating function
$D_2(x,q):=\sum_{n\geq 0} d_{n,2}(q)x^n$.
From the Theorem~\ref{degreegfg} we have the expression
$$D_2(x,q)=\frac{xq^4(2-x+2xq)}{1-x-q^2x^2}.$$
This generating function satisfies
$$D_2(x,q)-xD_2(x,q)-q^2x^2D_2(x,q)=2 q^4 x - q^4 x^2 + 2 q^5 x^2.$$
By comparing the  coefficient  of $x^n$, for  $n\geq 2$,  in the above equality, we obtain the recurrence relation for the sequence $d_{n,2}(q)$. On the other hand, the combinatorial sum \eqref{vertdeg2} is equivalent to 
\begin{align}\label{closedf2}
\sum_{i\geq 0} \left(2\binom{n-i-1}{i-1}q^{2i+3}+ \binom{n-i}{i-1}q^{2i+2} + \binom{n-i-1}{i-2}q^{2i+2} \right).
\end{align}
Denote the above summand  by $F(n,i)$, that is
\begin{align*}
F(n,i):=& 2\binom{n-i-1}{i-1}q^{2i+3}+ \binom{n-i}{i-1}q^{2i+2} + \binom{n-i-1}{i-2}q^{2i+2}\\
=&\frac{(n-i-1)!\left(n(2q +1) -4 i q  + 2q-1\right)}{(i-1)!(n-2i+1)!}q^{2i+2}.
\end{align*}
By the Zeilberger algorithm, $F(n,i)$ satisfies the relation
\begin{align}\label{rel2}
F(n+2,i)-F(n+1,i)-q^2F(n,i)=G(n,i+1)-G(n,i),
\end{align}
with the certificate 
$$R(n,i)=\frac{(n-i) (i-1)\left((4 i-6)q - n(2q+1)  + 1\right)}{(n+2-2 i) (n+3-2 i) \left(2 n q+n-4 i q+2 q-1\right)}.$$
If $f(n)$ denotes the combinatorial  sum in \eqref{closedf2}, then summing both sides of \eqref{rel2} with respect to $i$ yields $f(n+2)-f(n+1)-q^2f(n)=0$.  The sequences $f(n)$ and $d_{n,2}(q)$ satisfy the same recurrence relation and  have the same initial values, therefore these sequences coincides for all positive integers $n$.

\end{proof}

Notice that the  total vertices of degree 2 of all Fibonacci graphs in $\G_{n,2}$  is given by
$$\sum_{i=1}^{n}(i+3) \left(\binom{n-1-\lfloor \frac{i}{2}\rfloor}{\lfloor\frac{i-1}{2}\rfloor}+  \binom{n-2-\lfloor \frac{i-1}{2}\rfloor}{\lfloor\frac{i-2}{2}\rfloor}\right).$$

From a similar argument as in Theorem~\ref{closeformula} we can prove the following two theorems.

\begin{theorem}
For all $n\geq 3$ we have
$$d_{n,3}(q)=q^2(d_{n-1,3}(q) + d_{n-2,3}(q)),$$
with the initial values $d_{1,3}(q)=q^2+1$ and $d_{2,3}(q)=3q^2$. Moreover,  for all $n\geq 1$ we have the combinatorial formula 
$$d_{n,3}(q)=\sum_{i=0}^{n}\left((1+q^2)q^{2(n-1-i)}\binom{n-i-1}{i} + (2q^2-q^4)q^{2(n-1-i)}\binom{n-i-2}{i}\right).$$
\end{theorem}

\begin{theorem}
For all $n\geq 3$ we have
$$d_{n,4}(q)=d_{n-1,4}(q) + q^2d_{n-2,4}(q),$$
with the initial values $d_{1,2}(q)=2q^4$ and $d_{2,2}(q)=q^4 + 2 q^4$. Moreover,  for all $n\geq 1$ we have the combinatorial formula 
$$d_{n,4}(q)=\sum_{i=0}^{n+2} \left(\binom{n-1-\lfloor \frac{i-2}{2}\rfloor}{\lfloor\frac{i-3}{2}\rfloor}+  \binom{n-2-\lfloor \frac{i-3}{2}\rfloor}{\lfloor\frac{i-4}{2}\rfloor}\right)q^{i-3}.$$
\end{theorem}

Let $d_n$ denote the total number of vertices of the Fibonacci polyominoes in $\G_{n,2}$. In the following theorems we study the proportion between the sequences $d_{n,i}:=d_{n,i}(1)$ ($i=2, 3, 4$) and $d_n$.   Before,  we need the following result. 

\begin{theorem}[Asymptotics of linear recurrences,~\cite{Flajolet2}]\label{FlaTheo}
Assume that a rational generating function $f(x)/g(x)$, with $f(x)$ and $g(x)$ relatively prime and $g(0)\neq 0$, has a unique pole $1/\beta$ of smallest modulus. Then, if the multiplicity of $1/\beta$ is $\nu$, we have
$$[x^n]\frac{f(x)}{g(x)}\sim \nu \frac{(-\beta)^\nu f(1/\beta)}{g^{(\nu)}(1/\beta)}\beta^n n^{\nu-1}.$$

\end{theorem}

\begin{theorem}\label{asymp1}
Among total degree of  vertices of all graphs in $\G_{n,2}$, the proportion of those that are of degree 2 is asymptotically
\[
\lim_{n\to\infty}\frac{d_{n,2}}{d_n} = \frac{7-\sqrt{5}}{22}\approx 0.21654236.
\]
\end{theorem}
\begin{proof}
The generating functions of the sequences $d_n$ and $d_{n,2}$ are rational, therefore we can use the asymptotic analysis for linear recurrences. First, note that the unique pole $1/\beta$ of the rational generating function 
$$V_2(x)=\sum_{n \geq 0}d_nx^n=\frac{2 x (5 + x - 2 x^2 - x^3)}{(1 - x - x^2)^2}$$ is $\alpha=\frac{-1+\sqrt{5}}{2}$, with multiplicity 2. Therefore
$d_n\sim \frac{2(3+2\sqrt{5})}{5}\left(\frac{1+\sqrt{5}}{2}\right)^nn$. Similarly, we have $d_{n,2}\sim(\frac{1+\sqrt{5}}{2})^{n+1}n$. From these expressions we obtain the desired result.
\end{proof}

\begin{theorem}
Among total degree of  vertices of all graphs in $\G_{n,2}$, the proportion of those that are of degree 3 is asymptotically
\[
\lim_{n\to\infty}\frac{d_{n,3}}{d_n} = \frac{4+\sqrt{5}}{11}\approx 0.56691527.
\]
\end{theorem}
\begin{proof}
The unique pole $1/\beta$ of the rational generating function 
$$\sum_{n \geq 0}d_{n,3}x^n=\frac{2 (x - x^2 + x^3 - 2 x^4 + x^6)}{(1 - 2 x + x^3)^2}$$ is $\alpha=\frac{-1+\sqrt{5}}{2}$, with multiplicity 2. Therefore
$d_{n,2}\sim \frac{2(2+\sqrt{5})}{5}(\frac{1+\sqrt{5}}{2})^nn$. Since \linebreak $d_n\sim \frac{2(3+2\sqrt{5})}{5}(\frac{1+\sqrt{5}}{2})^nn$  we obtain the desired result.
\end{proof}

\begin{theorem}
Among total degree of  vertices of all graphs in $\G_{n,2}$, the proportion of those that are of degree 4 is asymptotically
\[
\lim_{n\to\infty}\frac{d_{n,4}}{d_n} = \frac{7-\sqrt{5}}{22}\approx 0.216542364.
\]
\end{theorem}

\section{Number of Hamiltonian $k$-bonacci graphs}
A \emph{Hamiltonian cycle} is a cycle that visits each vertex exactly once.  Let $G$ be a $k$-bonacci graph.  We define  $\Ham(G)=1$ if $G$ has a Hamiltonian cycle, and 0 otherwise.  If $\Ham(G)=1$, we say that $G$ is a Hamiltonian $k$-bonacci graph. 
   Define the generating function
$$H_k(x,q):=\sum_{n\geq 1} x^n \sum_{G\in \G_{n,k}}q^{\Ham(G)}, \quad k\geq 2,$$
 where $x$ marks the length of the corresponding $k$-bonacci word, i.e., the number of vertices in the bottom row of a graph minus one.
Similarly, we have the generating functions
$$H_{k,j}(x,q):=\sum_{n\geq 1} x^n \sum_{G\in \G_{n,k}^{(j)}}q^{\Ham(G)},  \quad \text{for} \quad j=1, 2.$$

\begin{theorem}
For all $k\geq 2$ we have
$$H_k(x,q)=\frac{x ((1 - x) x (1 - x^{2 \lfloor (k-1)/2\rfloor}) - 
   q (1 + x) (1 - 2 x + x^{k+1}) (-2 + x + x^{2 \lfloor k/2\rfloor}))}{(1 - 2 x + x^{k+1}) (1 - x - 2 x^2 + x^3 + x^{2(1 +  \lfloor k/2 \rfloor})}.$$
\end{theorem}
\begin{proof}
From the decomposition given in Figure~\ref{deco} we have the functional equation
\begin{align*}
H_{k,1}(x,q)&=\underbrace{qx}_{(1)} + \underbrace{x(H_{k,1}(x,q)+H_{k,2}(x,q))}_{(2)+(3)} 
\end{align*}
The polyomino given in the decomposition $(4)$ corresponds to the grid graph $P_3 \times P_i$, for $1 < i \leq k$. It is known\footnote{~See, for example the discussion on ``Math Stackexchange": \url{https://math.stackexchange.com/questions/1699203/hamilton-paths-cycles-in-grid-graphs}} that the grid graph $P_n \times P_m$ has a hamiltonian cycle if and only if at least one of $m$ or $n$ is even
or $m=n=1$. Therefore, the graph $P_3 \times P_i$ has a hamiltonian cycle if and only if $i$ is even. In this case, the generating function is given by
$$T_k(x,q)=\underbrace{q\sum_{j=1}^{\lfloor k/2 \rfloor}x^{2j-1} + \sum_{j=1}^{\lfloor (k-1)/2 \rfloor}x^{2j}}_{(4)} =\frac{x\left(q + x - x^{2\lfloor\frac{k-1}{2} \rfloor+1} -  qx^{2\lfloor\frac{k}{2} \rfloor} \right)}{1-x^2}.
$$
Therefore, we obtain the functional equation
\begin{align*}
H_{k,2}(x,q)&=\underbrace{T_k(x,q)}_{(4)} +\underbrace{\left(\sum_{j=1}^{\lfloor k/2 \rfloor}x^{2j-1}\right) H_{k,1}(x,q) +\left(\sum_{j=1}^{\lfloor(k-1)/2 \rfloor}x^{2j}\right) H_{k,1}(x,1)}_{(5)}\\
&=T_k(x,q) + \frac{x(1-x^{2\lfloor k/2 \rfloor})}{1-x^2}H_{k,1}(x,q) + \frac{x^2(1-x^{2\lfloor (k-1)/2 \rfloor})}{1-x^2}\frac{x(1-x^k)}{1-2x+x^{k+1}}.
\end{align*}
Notice that $H_{k,1}(x,1)$ is the generating function of the total number of $k$-bonacci words end in \texttt{1}.
Solving the above system of equation  we obtain the desired result.
\end{proof}

\begin{corollary}
The generating function for the total number of Hamiltonian $k$-bonacci graphs is
$$H_{k}(x):=\frac{x(1+x)(2-x-x^{2\lfloor k/2 \rfloor})}{1 - x - 2 x^2 + x^3 + x^{2 \lfloor k/2\rfloor +2}}.$$
\end{corollary}

In particular, for $2 \le k \le 7$
we obtain the following generating functions:
\begin{align*}
H_{2}(x)&=H_{3}(x)=\frac{x (2 + x)}{1 -  x -x^2},\\
H_{4}(x)&=H_{5}(x)=\frac{x (2 + x + x^2 + x^3)}{1 - x - x^2 - x^4},\\
H_{6}(x)&=H_{7}(x)=\frac{x (2 + x + x^2 + x^3 + x^4 + x^5)}{1 - x - x^2 - x^4 - x^5}.
   \end{align*}

Note that every Fibonacci graph admits a Hamiltonian cycle,
just walk along its border.
The grid $P_3 \times P_3$ does not have a Hamiltonian cycle, so
Hamiltonian 3-bonacci graphs cannot have grids $P_3 \times P_3$ as induced subgraphs. Thus, they are precisely 2-bonacci.
In general, Hamiltonian $2k$-bonacci equals
Hamiltonian $(2k+1)$-bonacci graphs, because any Hamiltonian $(2k+1)$-bonacci graph cannot
contain $P_{2k+1} \times P_3$ as
induced subgraph, and there are no $1^{2k}$ factors in the corresponding binary words.

\medskip

\noindent MSC2020: 11B39, 05A15, 05A19. 


\begin{thebibliography}{99}

\bibitem{Allou} J.P.~Allouche and J.~Johnson, Narayana's cows and delayed morphisms,  In: \emph{Articles of 3rd Computer Music Conference JIM96}, France, (1996).

\bibitem{Baril} J.-L.~Baril, S.~Kirgizov, and V.~Vajnovszki, Gray codes for Fibonacci $q$-decreasing words, \emph{Theoret. Comput. Sci.}
\textbf{927} (2022), 120--132.

\bibitem{NIV} D.~Beauquier and M.~Nivat, On translating one polyomino to tile the plane, \emph{Discrete Comput. Geom.} \textbf{6} (1991), 575--592.

\bibitem{Bernini} A.~Bernini, Restricted binary strings and generalized Fibonacci numbers, In: \emph{International Workshop on Cellular Automata and Discrete Complex Systems}, Springer (2017),  32--43. 

\bibitem{BLE3} A.~Blecher, C.~Brennan, and A.~Knopfmacher, Combinatorial parameters in bargraphs, \emph{Quaest. Math.} \textbf{39} (2016), 619--635.

\bibitem{CallManRam} D.~Callan, T.~Mansour, and J.L.~Ram\'irez,  Statistics on bargraphs of Catalan words, \emph{J. Autom. Lang. Comb.}  \textbf{26} (2021), 177--196.

\bibitem{BLO2} A.~Blondin-Massé, S.~Brlek, A.~Garon, and S.~Labbé, Two infinite families of polyominoes that tile the plane by translation in two distinct ways, \emph{Theoret. Comput. Sci.} \textbf{412} (2011), 4778--4786.

\bibitem{RAM} R.~De~Castro, A.~Ram\'{\i}rez, and J.~L.~Ram\'{\i}rez, Applications in enumerative combinatorics of infinite weighted
automata and graphs, \emph{Sci. Ann. Comput. Sci.} \textbf{24} (2014), 137--171.

\bibitem{DEL} M.~Delest, Algebraic languages: a bridge between combinatorics and computer science, \emph{Discrete Math. Theor. Comput. Sci.} \textbf{24} (1996), 71--88.

\bibitem{OME1}  \" O. E\u{g}ecio\u{g}lu, Statistics on restricted Fibonacci words, \emph{Trans. Combin.} \textbf{10}(1) (2020), 31--42.

\bibitem{OME2} \" O. E\u{g}ecio\u{g}lu and V. Ir\v{s}i\v{c}, Fibonacci-run graphs I: Basic properties, \emph{Discrete Applied Math.} \textbf{295} (2021), 70--84.


\bibitem{OME3} \" O. E\u{g}ecio\u{g}lu and V. Ir\v{s}i\v{c}, Fibonacci-run graphs II: Degree sequences, \emph{Discrete Applied Math.} \textbf{300} (2021), 56--71.


\bibitem{flajolet} P.~Flajolet and R.~Sedgewick,  \emph{Analytic Combinatorics}, Cambridge University Press, Cambridge, 2009.

 \bibitem{Book1} A.~J.~Guttmann (Ed.), \emph{Polygons, Polyominoes and Polycubes}, Lecture Notes in Physics 775. Springer, Heidelberg, Germany, 2009.
 
\bibitem{Hsu} W.~J.~Hsu, Fibonacci cubes - new interconnection topology, \emph{Parallel and Distributed Systems,
IEEE Transactions}  \textbf{4} (1993), 3--12.


\bibitem{Klav2} S.~Klav\v zar, Structures of Fibonacci cubes: a survey, \emph{J. Comb. Optim.} \textbf{25} (2013), 505--522.



\bibitem{Klav5} S.~Klav\v zar, M.~Mollard, and M.~Petkov\v{s}ek, The degree sequence of Fibonacci and Lucas cubes, \emph{Discrete Math.} \textbf{311} (14) (2011), 1310--1322.

\bibitem{knuth3}
  D.~E.~Knuth, \emph{The Art of Computer Programming, Volume 3: Sorting and Searching}, 2nd ed.  Addison-Wesley, 1998.

\bibitem{Koshy} T.~Koshy, \emph{Fibonacci and Lucas  Number with Applications},  John Wiley \& Sons, 2001.  

  \bibitem{ManSha} T.~Mansour and A.~Sh.~Shabani, Enumerations on bargraphs, \emph{Discrete Math. Lett.} \textbf{2} (2019), 65--94.

\bibitem{ABZ} M.~Petko\v{v}sek, H.~Wilf, and D.~Zeilberger, \emph{A=B}, A. K. Peters, Ltd. 1996.
  
  \bibitem{RAMRUB}
   J.~L.~Ram\'irez, G.~N.~Rubiano, and R.~De~Castro, A generalization of the Fibonacci word Fractal and the Fibonacci snowflake, \emph{Theoret. Comput. Sci.} \textbf{528}(2014), 40--56.
  
\bibitem{Flajolet2} R.~Sedgewick and P.~Flajolet,  \emph{An Introduction to the Analysis of Algorithms} 2nd ed., Addison-Wesley, 2013.
 
\bibitem{OEIS}
OEIS Foundation Inc.~(2022),
{\em The On-Line Encyclopedia of Integer Sequences},
\url{https://oeis.org}.


\bibitem{Vaj} V.~Vajnovszki, A loopless generation of bitstrings without $p$ consecutive ones,  In:  \emph{Combinatorics, Computability and Logic}. \emph{Discrete Math. Theor. Comput. Sci.}  Springer (2001), 227--240.



\end{thebibliography}
\end{document}